\documentclass[12pt]{article}%
\usepackage{amsmath}
\usepackage{amsfonts}
\usepackage{amssymb}
\usepackage{graphicx}%
\newtheorem{theorem}{Theorem}

\newtheorem{corollary}[theorem]{Corollary}

\begin{document}

\title{Asymptotic analysis of the Hermite polynomials from their
differential-difference equation}
\author{Diego Dominici \thanks{e-mail: dominicd@newpaltz.edu}\\Department of Mathematics\\State University of New York at New Paltz\\75 S. Manheim Blvd. Suite 9\\New Paltz, NY 12561-2443\\USA\\Phone: (845) 257-2607\\Fax: (845) 257-3571 }
\maketitle

\begin{abstract}
We analyze the Hermite polynomials $H_{n}(x)$ and their zeros asymptotically,
as $n\rightarrow\infty.$ We obtain asymptotic approximations from the
differential-difference equation which they satisfy, using the ray method. We
give numerical examples showing the accuracy of our formulas.

\end{abstract}


MSC-class: 33C45 (Primary) 34E05, 34E20 (Secondary)

\section{Introduction}

It would be difficult to find a more ubiquitous polynomial family than the
Hermite polynomials $H_{n}(x),$ defined by the Rodrigues formula%
\[
H_{n}(x)=\left(  -1\right)  ^{n}\exp\left(  x^{2}\right)  \frac{d^{n}}{dx^{n}%
}\exp\left(  -x^{2}\right)  ,\quad n=0,1,2,\ldots.
\]
They appear in several problems of mathematical physics \cite{MR0350075}, the
most important probably being the solution of the Schr\"{o}dinger equation
\cite{MR1270235}, \cite{MR1465675}. Being the limiting case of several
families of classical orthogonal polynomials \cite{koekoek94askeyscheme}, they
are of fundamental importance in asymptotic analysis \cite{MR1723073},
\cite{MR1805994}.

The Hermite polynomials satisfy the orthogonality condition%
\[%
{\displaystyle\int\limits_{-\infty}^{\infty}}
e^{-x^{2}}H_{m}(x)H_{n}(x)dx=\sqrt{\pi}2^{n}n!\delta_{mn,}%
\]
the differential-difference equation%
\begin{equation}
H_{n+1}+H_{n}^{\prime}=2xH_{n}, \label{diffdiff}%
\end{equation}
and the reflection formula
\begin{equation}
H_{n}(-x)=\left(  -1\right)  ^{n}H_{n}(x). \label{reflex}%
\end{equation}

The zeros of the Hermite polynomials have several applications, notably in
Gauss' quadrature formula for numerical integration \cite{MR0026415},
\cite{MR0036084}, \cite{MR0048901}. Several properties and their asymptotic
behavior were studied in \cite{MR502800}, \cite{MR2059744}, \cite{MR0457815},
\cite{MR522774}, \cite{MR1342385}, \cite{MR1343543}, \cite{MR1343062},
\cite{MR1303614} and \cite{MR1329018}.

The asymptotic behavior of $H_{n}(x)$ was studied by M. Plancherel and W.
Rotach in \cite{MR1509395} using the method that now bears their name. F. W.
J. Olver \cite{MR0109898} obtained asymptotic expansions for the Hermite
polynomials as a consequence of his WKB analysis of the differential equation
satisfied by the Parabolic Cylinder function $D_{\nu}\left(  z\right)  ,$
related to $H_{n}(x)$ by
\[
H_{n}(x)=2^{\frac{n}{2}}\exp\left(  \frac{x^{2}}{2}\right)  D_{n}\left(
\sqrt{2}x\right)  .
\]
A similar analysis using perturbation techniques was carried on by A. Voznyuk
in \cite{MR924562}.

As an application of the results from his doctoral thesis on the
multiplication-interpolation method, L. Heflinger \cite{MR0086169} established
asymptotic series for the Hermite polynomials. In \cite{MR0146423}, M. Wyman
derived asymptotic formulas for $H_{n}(x)$ based on one of their integral representations.

In this paper we shall take a different approach and analyze the
differential-difference equations that the Hermite polynomials satisfy
(\ref{diffdiff}) using the techniques presented in \cite{MR0361328}. A similar
method (which we may call the discrete WKB method) has been applied to the
solution of difference equations \cite{MR1373150}, \cite{MR0225511},
\cite{MR1164992}, \cite{MR804827} and it is currently being extended
\cite{MR2103370}, \cite{MR1896091}, \cite{MR1971216}, \cite{MR2114641}, to
include difference equations with turning points. Another type of analysis,
based on perturbation techniques, was considered by C. Lange and R. Miura in
\cite{MR659409}, \cite{MR804001}, \cite{MR804002}, \cite{MR1093425},
\cite{MR1260916}, and \cite{MR1260917}.

\section{Asymptotic analysis}

We consider the approximation
\begin{equation}
H_{n}(x)\sim\exp\left[  f(x,n)+g(x,n)\right]  ,\quad n\rightarrow
\infty\label{asympH}%
\end{equation}
where%
\begin{equation}
g=o(f),\quad n\rightarrow\infty. \label{gof}%
\end{equation}
Note that since $H_{0}(x)=1,$ we must have%
\begin{equation}
f(x,0)=0 \label{f(0)}%
\end{equation}
and%
\begin{equation}
g(x,0)=0. \label{g(0)}%
\end{equation}
Using (\ref{asympH}) in (\ref{diffdiff}), we have%
\begin{gather}
\exp\left(  f+\frac{\partial f}{\partial n}+\frac{1}{2}\frac{\partial^{2}%
f}{\partial n^{2}}+g+\frac{\partial g}{\partial n}\right) \label{asymptH1}\\
+\left(  \frac{\partial f}{\partial x}+\frac{\partial g}{\partial x}\right)
\exp\left(  f+g\right)  =2x\exp\left(  f+g\right)  ,\nonumber
\end{gather}
where we have used
\[
f(x,n+1)=f(x,n)+\frac{\partial f}{\partial n}(x,n)+\frac{1}{2}\frac
{\partial^{2}f}{\partial n^{2}}(x,n)+\cdots.
\]
Simplifying (\ref{asymptH1}) and taking (\ref{gof}) into account we obtain, to
leading order,%
\begin{equation}
\exp\left(  \frac{\partial f}{\partial n}\right)  +\frac{\partial f}{\partial
x}=2x. \label{eqf}%
\end{equation}
Using (\ref{eqf}) in (\ref{asymptH1}) we get%
\[
\exp\left(  \frac{1}{2}\frac{\partial^{2}f}{\partial n^{2}}+\frac{\partial
g}{\partial n}\right)  +\frac{\partial g}{\partial x}\exp\left(
-\frac{\partial f}{\partial n}\right)  =1,
\]
or, to leading order,%
\begin{equation}
\frac{1}{2}\frac{\partial^{2}f}{\partial n^{2}}+\frac{\partial g}{\partial
n}+\frac{\partial g}{\partial x}\exp\left(  -\frac{\partial f}{\partial
n}\right)  =0. \label{eqg}%
\end{equation}

\subsection{The ray expansion}

To solve (\ref{eqf}) we use the method of characteristics, which we briefly
review. Given the first order partial differential equation%
\[
F\left(  x,n,f,p,q\right)  =0,
\]
where
\[
\ p=\frac{\partial f}{\partial x},\quad q=\frac{\partial f}{\partial n},
\]
we search for a solution \ $f(x,n)$ by solving the system of \textquotedblleft
characteristic equations\textquotedblright\
\begin{align*}
\frac{dx}{dt}  &  =\frac{\partial F}{\partial p},\quad\frac{dn}{dt}%
=\frac{\partial F}{\partial q},\\
\frac{dp}{dt}  &  =-\frac{\partial F}{\partial x}-p\frac{\partial F}{\partial
f},\quad\frac{dq}{dt}=-\frac{\partial F}{\partial n}-q\frac{\partial
F}{\partial f},\\
\frac{df}{dt}  &  =p\frac{\partial F}{\partial p}+q\frac{\partial F}{\partial
q},
\end{align*}
where we now consider $\left\{  x,n,f,p,q\right\}  $ to all be functions of
the variables $t$ and $s.$

For (\ref{eqf}), we have%
\begin{equation}
F\left(  x,n,f,p,q\right)  =e^{q}+p-2x \label{F}%
\end{equation}
and therefore the characteristic equations are%
\begin{equation}
\frac{dx}{dt}=1,\quad\frac{dn}{dt}=e^{q},\quad\frac{dp}{dt}=2,\quad\frac
{dq}{dt}=0, \label{sysray}%
\end{equation}
and
\begin{equation}
\frac{df}{dt}=p+qe^{q}. \label{eqfray}%
\end{equation}
Solving (\ref{sysray}) subject to the initial conditions%
\begin{equation}
x(0,s)=s,\quad n(0,s)=0,\quad q(0,s)=A(s), \label{IC}%
\end{equation}
we obtain%
\begin{equation}
x=t+s,\quad n=te^{A},\quad p=2t+2s-e^{A},\quad q=A, \label{sol1}%
\end{equation}
where we have used%
\[
0=\left.  F\right\vert _{t=0}=e^{A}+p(0,s)-2s.
\]
From (\ref{f(0)}) and (\ref{IC}) we have%
\begin{equation}
f(0,s)=0, \label{ft=0}%
\end{equation}
\ which implies%
\begin{align*}
0  &  =\frac{d}{ds}f(0,s)=\left[  \frac{\partial f}{\partial x}\frac{\partial
x}{\partial s}+\frac{\partial f}{\partial n}\frac{\partial n}{\partial
s}\right]  _{t=0}\\
&  =p(0,s)\times1+q(0,s)\times0=2s-e^{A}.
\end{align*}
Thus, \ $A(s)=\ln(2s)$ \ and (\ref{sol1}) becomes%
\begin{equation}
x=t+s,\quad n=2ts,\quad p=2t,\quad q=\ln(2s), \label{sol2}%
\end{equation}
with $t\geq0$ $\ $and $\ s>0.$ Since $s>0,$ we shall consider only the region
$x>0$ for now. Using (\ref{sol2}) in (\ref{eqfray}) and taking (\ref{ft=0})
into account, we obtain%
\begin{equation}
f(t,s)=t^{2}+2s\ln(2s)t. \label{f(t,s)}%
\end{equation}

Solving for $t$ and $s$ in terms of $x$ and $n$ in (\ref{sol2}), we get%
\begin{equation}
t=\frac{x}{2}\pm\frac{1}{2}\sigma,\quad s=\frac{x}{2}\mp\frac{1}{2}%
\sigma\label{ts}%
\end{equation}
with%
\begin{equation}
\sigma=\sqrt{x^{2}-2n}. \label{sigma}%
\end{equation}
For $\sigma$ to be a real number, we shall impose the condition $x>\sqrt{2n.}$
Since (for a fixed value of $n)$ we have $t\rightarrow0$ as $x\rightarrow
\infty,$ we consider the solution%
\begin{equation}
t=\frac{x}{2}-\frac{1}{2}\sigma,\quad s=\frac{x}{2}+\frac{1}{2}\sigma.
\label{ts1}%
\end{equation}
Replacing (\ref{ts1}) in (\ref{f(t,s)}) we obtain%
\begin{equation}
f(x,n)=\frac{x^{2}-\sigma x-n}{2}+n\ln\left(  x+\sigma\right)  ,\quad
x>\sqrt{2n}. \label{f}%
\end{equation}

We shall now find $g(x,n)$. Using (\ref{f}) in (\ref{eqg}), we get%
\[
-\frac{1}{2\sigma\left(  x+\sigma\right)  }+\frac{\partial g}{\partial
n}+\frac{\partial g}{\partial x}\frac{1}{\left(  x+\sigma\right)  }=0,
\]
or
\begin{equation}
\left(  x+\sigma\right)  \frac{\partial g}{\partial n}+\frac{\partial
g}{\partial x}=\frac{1}{2\sigma}. \label{dg}%
\end{equation}
Solving (\ref{dg}), we obtain%
\[
g(x,n)=\frac{1}{2}\ln\left(  -2\frac{x^{2}-n+x\sigma}{x^{2}-2n+x\sigma
}\right)  +C\left(  x+\sigma\right)  ,
\]
where $C(x)$ is a function to be determined. Imposing the condition
(\ref{g(0)}), we have%
\[
0=g(x,0)=\frac{1}{2}\ln(-2)+C(2x).
\]
Thus,
\begin{equation}
g(x,n)=\frac{1}{2}\ln\left[  \frac{1}{2}\left(  \frac{x}{\sigma}+1\right)
\right]  . \label{g}%
\end{equation}

We summarize our results in the following theorem.

\begin{theorem}
In the region $x>\sqrt{2n},$ the Hermite polynomials admit the asymptotic
representation%
\begin{align}
H_{n}(x)\sim\Phi_{1}(x,n)=  &  \exp\left[  \frac{x^{2}-\sigma x-n}{2}%
+n\ln\left(  \sigma+x\right)  \right] \label{Phi1}\\
&  \times\sqrt{\frac{1}{2}\left(  1+\frac{x}{\sigma}\right)  },\quad
n\rightarrow\infty,\nonumber
\end{align}
where $\sigma(x,n)$ was defined in (\ref{sigma}).
\end{theorem}

Using the reflection formula (\ref{reflex}) we can extend our result to the
region $-x>\sqrt{2n}$ and obtain:

\begin{corollary}
In the region $x<-\sqrt{2n},$ the Hermite polynomials admit the asymptotic
representation%
\begin{align}
H_{n}(x)\sim\Phi_{2}(x,n)=  &  \left(  -1\right)  ^{n}\exp\left[  \frac
{x^{2}+\sigma x-n}{2}+n\ln\left(  \sigma-x\right)  \right] \label{Phi2}\\
&  \times\sqrt{\frac{1}{2}\left(  1-\frac{x}{\sigma}\right)  },\quad
n\rightarrow\infty.\nonumber
\end{align}

\end{corollary}

\begin{figure}[ptb]
\begin{center}
\rotatebox{270} {\resizebox{!}{5in}{\includegraphics{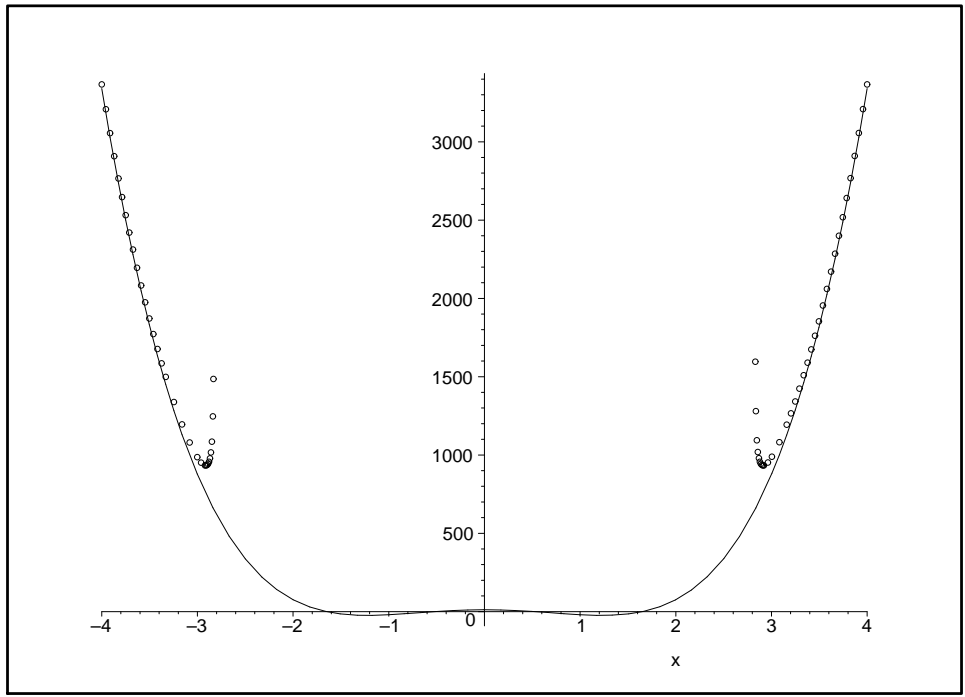}}}
\end{center}
\caption{A comparison of $H_{4}(x)$ (solid curve) and the asymptotic
approximations $\Phi_{1}(x,4)$ and $\Phi_{2}(x,4)$ (ooo).}%
\label{fexp}%
\end{figure}

To illustrate the accuracy of our results, in Figure \ref{fexp} we graph
$H_{4}(x)$ and the asymptotic approximations $\Phi_{1}(x,4)$ and $\Phi
_{2}(x,4).$

\subsection{The transition layer}

We shall now find an asymptotic approximation for $\left\vert x\right\vert
\approx\sqrt{2n}.$ We will consider the case $x\approx\sqrt{2n}$ and find the
corresponding result for $x\approx-\sqrt{2n}$ by using (\ref{reflex}). From
(\ref{Phi1}) we have%
\[
\Phi_{1}(x,n)\sim\exp\left[  \frac{n}{2}\ln\left(  2n\right)  -\frac{3}%
{2}n+\sqrt{2n}x\right]  ,\quad x\longrightarrow\sqrt{2n}^{+}.
\]
We define the function $G_{n}(x)$ by%
\begin{equation}
H_{n}(x)=\exp\left[  \frac{n}{2}\ln\left(  2n\right)  -\frac{3}{2}n+\sqrt
{2n}x\right]  G_{n}(x). \label{G}%
\end{equation}
Using (\ref{G}) in (\ref{diffdiff}) we get%
\begin{gather}
\exp\left[  \frac{n+1}{2}\ln\left(  2n+2\right)  -\frac{n}{2}\ln\left(
2n\right)  -\frac{3}{2}+\left(  \sqrt{2\left(  n+1\right)  }-\sqrt{2n}\right)
x\right]  G_{n+1}\label{G1}\\
+\sqrt{2n}G_{n}+G_{n}^{\prime}=2xG_{n}(x).\nonumber
\end{gather}

We introduce the stretch variable $\beta>0$ defined by%
\begin{equation}
x=\sqrt{2n}+\frac{\beta}{n^{\frac{1}{6}}}\label{beta}%
\end{equation}
and the function $\Lambda\left(  \beta\right)  $ defined by%
\begin{equation}
G_{n}(x)=\Lambda\left[  \left(  x-\sqrt{2n}\right)  n^{\frac{1}{6}}\right]
.\label{lambda}%
\end{equation}
From (\ref{beta}) we have%
\begin{gather}
\exp\left[  \frac{n+1}{2}\ln\left(  2n+2\right)  -\frac{n}{2}\ln\left(
2n\right)  -\frac{3}{2}+\left(  \sqrt{2\left(  n+1\right)  }-\sqrt{2n}\right)
x\right]  \label{exp}\\
\sim\sqrt{2n}+\beta n^{-\frac{1}{6}},\quad n\rightarrow\infty.\nonumber
\end{gather}
Using (\ref{beta}) in (\ref{lambda}) we obtain%
\begin{align}
G_{n+1}(x) &  =\Lambda\left[  \left(  \sqrt{2n}-\sqrt{2n+1}+\frac{\beta
}{n^{\frac{1}{6}}}\right)  \left(  n+1\right)  ^{\frac{1}{6}}\right]
\label{G2}\\
&  \sim\Lambda\left(  \beta\right)  -\frac{1}{\sqrt{2}}\Lambda^{\prime}\left(
\beta\right)  n^{-\frac{1}{3}}+\frac{1}{4}\Lambda^{\prime\prime}\left(
\beta\right)  n^{-\frac{2}{3}},\quad n\rightarrow\infty,\nonumber
\end{align}
and%
\begin{equation}
\sqrt{2n}G_{n}+G_{n}^{\prime}-2xG_{n}(x)=-\sqrt{2n}\Lambda\left(
\beta\right)  +\Lambda^{\prime}\left(  \beta\right)  n^{\frac{1}{6}}%
-2\beta\Lambda\left(  \beta\right)  n^{-\frac{1}{6}}.\label{G3}%
\end{equation}
Using (\ref{exp}), (\ref{G2}) and (\ref{G3}) in (\ref{G1}) we obtain, to
leading order, the Airy equation%
\begin{equation}
\Lambda^{\prime\prime}\left(  \beta\right)  =2\sqrt{2}\beta\Lambda\left(
\beta\right)  .\label{airy}%
\end{equation}
Thus,
\begin{equation}
\Lambda\left(  \beta\right)  =C_{1}\operatorname{Ai}\left(  \sqrt{2}%
\beta\right)  +C_{2}\operatorname{Bi}\left(  \sqrt{2}\beta\right)
,\label{lambda1}%
\end{equation}
where $\operatorname{Ai}\left(  \cdot\right)  $ and $\operatorname{Bi}\left(
\cdot\right)  $ denote the Airy functions and $C_{1},$ $C_{2}$ are to be
determined. Replacing (\ref{beta}) and (\ref{lambda1}) in (\ref{G}) we have%
\begin{equation}
H_{n}(x)\sim\exp\left[  \frac{n}{2}\ln\left(  2ne\right)  +\sqrt{2}\beta
n^{\frac{1}{3}}\right]  \left[  C_{1}\operatorname{Ai}\left(  \sqrt{2}%
\beta\right)  +C_{2}\operatorname{Bi}\left(  \sqrt{2}\beta\right)  \right]
.\label{Htran}%
\end{equation}

To find $C_{1},$ $C_{2}$ we shall match (\ref{Htran}) with (\ref{Phi1}). Using
(\ref{beta}) in (\ref{Phi1}) we get%
\begin{equation}
\Phi_{1}(x,n)\sim\exp\left[  \frac{n}{2}\ln\left(  2ne\right)  +\sqrt{2}\beta
n^{\frac{1}{3}}-\frac{2^{\frac{7}{4}}}{3}\beta^{\frac{3}{2}}\right]
2^{-\frac{5}{8}}\beta^{-\frac{1}{4}}n^{\frac{1}{6}}, \label{phi10}%
\end{equation}
as $\beta\rightarrow0.$ Using (\ref{beta}) and the well known asymptotic
expansions of the Airy functions%
\begin{align*}
\operatorname{Ai}\left(  x\right)   &  \sim\frac{1}{2\sqrt{\pi}}\exp\left(
-\frac{2}{3}x^{\frac{3}{2}}\right)  x^{-\frac{1}{4}},\quad x\rightarrow
\infty\\
\operatorname{Bi}\left(  x\right)   &  \sim\frac{1}{\sqrt{\pi}}\exp\left(
\frac{2}{3}x^{\frac{3}{2}}\right)  x^{-\frac{1}{4}},\quad x\rightarrow\infty,
\end{align*}
in (\ref{Htran}) we have%
\begin{equation}
\exp\left[  \frac{n}{2}\ln\left(  2n\right)  -\frac{3}{2}n+\sqrt{2n}x\right]
=\exp\left\{  \frac{n}{2}\left[  1+\ln\left(  2n\right)  \right]  +\sqrt
{2}\beta n^{\frac{1}{3}}\right\}  \label{exp1}%
\end{equation}
and%
\begin{gather}
C_{1}\operatorname{Ai}\left(  \sqrt{2}\beta\right)  +C_{2}\operatorname{Bi}%
\left(  \sqrt{2}\beta\right)  \sim\label{ai1}\\
\frac{C_{1}}{\sqrt{\pi}2^{\frac{9}{8}}\beta^{\frac{1}{4}}}\exp\left(
-\frac{2^{\frac{7}{4}}}{3}\beta^{\frac{3}{2}}\right)  +\frac{C_{2}}{\sqrt{\pi
}2^{\frac{1}{8}}\beta^{\frac{1}{4}}}\exp\left(  \frac{2^{\frac{7}{4}}}{3}%
\beta^{\frac{3}{2}}\right)  ,\quad\beta\rightarrow\infty.\nonumber
\end{gather}
Matching (\ref{phi10}) to (\ref{exp1}) and (\ref{ai1}), we conclude that%
\begin{equation}
C_{1}=\sqrt{2\pi}n^{\frac{1}{6}},\quad C_{2}=0. \label{C1C2}%
\end{equation}
This completes the analysis. Combining the results above, we have the
following result:

\begin{theorem}
For $x\approx\sqrt{2n},$ the Hermite polynomials have the asymptotic
representation%
\begin{align}
H_{n}(x)  &  \sim\Phi_{3}(x,n)=\exp\left[  \frac{n}{2}\ln\left(  2n\right)
-\frac{3}{2}n+\sqrt{2n}x\right] \label{phi3}\\
&  \times\sqrt{2\pi}n^{\frac{1}{6}}\operatorname{Ai}\left[  \sqrt{2}\left(
x-\sqrt{2n}\right)  n^{\frac{1}{6}}\right]  ,\quad n\rightarrow\infty
.\nonumber
\end{align}

\end{theorem}

Use of the reflection formula (\ref{reflex}) provides the corresponding result
for $x\approx-\sqrt{2n}.$

\begin{corollary}
For $x\approx-\sqrt{2n},$ the Hermite polynomials have the asymptotic
representation%
\begin{align}
H_{n}(x)  &  \sim\Phi_{4}(x,n)=\left(  -1\right)  ^{n}\exp\left[  \frac{n}%
{2}\ln\left(  2n\right)  -\frac{3}{2}n-\sqrt{2n}x\right] \label{phi4}\\
&  \times\sqrt{2\pi}n^{\frac{1}{6}}\operatorname{Ai}\left[  -\sqrt{2}\left(
x+\sqrt{2n}\right)  n^{\frac{1}{6}}\right]  ,\quad n\rightarrow\infty
.\nonumber
\end{align}

\end{corollary}

\subsection{The oscillatory region}

We now study the region bounded by the curve $n=\frac{x^{2}}{2},$ where the
zeros of $H_{n}(x)$ are located. In this region, the solution is a linear
combination of (\ref{Phi1}) and (\ref{Phi2})%
\[
H_{n}(x)\sim\Phi_{5}(x,n)\equiv K_{1}\Phi_{1}(x,n)+K_{2}\Phi_{2}(x,n),\quad
n\rightarrow\infty
\]
with $\left\vert x\right\vert <\sqrt{2n}$ and $K_{1}$,$K_{2}$ are constants to
be determined. We shall require $\Phi_{5}(x,n)$ to match $\Phi_{3}(x,n)$
asymptotically in the local variable $\beta$, i.e., it must satisfy the
limiting condition%
\[
\underset{\beta\rightarrow0}{\lim}\Phi_{5}(\beta,n)=\underset{\beta
\rightarrow-\infty}{\lim}\Phi_{3}(\beta,n).
\]
\ 

Writing (\ref{phi3}) in terms of $\beta,$ we have%
\begin{equation}
\Phi_{3}(\beta,n)=\exp\left[  \frac{n}{2}\ln\left(  2ne\right)  +\sqrt{2}\beta
n^{\frac{1}{3}}\right]  \sqrt{2\pi}n^{\frac{1}{6}}\operatorname{Ai}\left(
\sqrt{2}\beta\right)  . \label{phi3beta}%
\end{equation}
Using the asymptotic formula%
\[
\operatorname{Ai}\left(  x\right)  \sim\frac{1}{\sqrt{\pi}}\sin\left[
\frac{2}{3}\left(  -x\right)  ^{\frac{3}{2}}+\frac{\pi}{4}\right]  \left(
-x\right)  ^{-\frac{1}{4}},\quad x\rightarrow-\infty
\]
in (\ref{phi3beta}) we get%
\begin{gather*}
\Phi_{3}(\beta,n)\sim\exp\left[  \frac{n}{2}\ln\left(  2ne\right)  +\sqrt
{2}\beta n^{\frac{1}{3}}\right]  2^{\frac{3}{8}}n^{\frac{1}{6}}\\
\times\sin\left[  \frac{1}{3}2^{\frac{7}{4}}\left(  -\beta\right)  ^{\frac
{3}{2}}+\frac{\pi}{4}\right]  \left(  -\beta\right)  ^{-\frac{1}{4}}%
,\quad\beta\rightarrow-\infty,
\end{gather*}
which can be rewritten as%
\begin{gather}
\Phi_{3}(\beta,n)\sim2^{-\frac{5}{8}}n^{\frac{1}{6}}\beta^{-\frac{1}{4}}%
\exp\left[  \frac{n}{2}\ln\left(  2ne\right)  +\sqrt{2}\beta n^{\frac{1}{3}%
}\right] \label{phi31}\\
\times\left[  \exp\left(  -\frac{1}{3}2^{\frac{7}{4}}\beta^{\frac{3}{2}%
}\right)  +\mathrm{i}\exp\left(  \frac{1}{3}2^{\frac{7}{4}}\beta^{\frac{3}{2}%
}\right)  \right]  ,\quad\beta\rightarrow-\infty.\nonumber
\end{gather}

Using (\ref{beta}) in (\ref{Phi1}), we have%
\begin{equation}
\Phi_{1}(\beta,n)\sim\exp\left[  \frac{n}{2}\ln\left(  2ne\right)  +\sqrt
{2}\beta n^{\frac{1}{3}}-\frac{1}{3}2^{\frac{7}{4}}\beta^{\frac{3}{2}}\right]
2^{-\frac{5}{8}}n^{\frac{1}{6}}\beta^{-\frac{1}{4}},\quad\beta\rightarrow0.
\label{phi11}%
\end{equation}
Similarly, using (\ref{beta}) in (\ref{Phi2}), we obtain%
\begin{equation}
\Phi_{2}(\beta,n)\sim\exp\left[  \frac{n}{2}\ln\left(  2ne\right)  +\sqrt
{2}\beta n^{\frac{1}{3}}+\frac{1}{3}2^{\frac{7}{4}}\beta^{\frac{3}{2}}\right]
2^{-\frac{5}{8}}n^{\frac{1}{6}}\beta^{-\frac{1}{4}}\mathrm{i},\quad
\beta\rightarrow0, \label{phi21}%
\end{equation}
where we have used
\[
\left(  -1\right)  ^{n}\exp\left[  \frac{x^{2}+\sigma x-n}{2}+n\ln\left(
\sigma-x\right)  \right]  =\exp\left[  \frac{x^{2}+\sigma x-n}{2}+n\ln\left(
x-\sigma\right)  \right]  .
\]
Comparing (\ref{phi31}) with (\ref{phi11}) and (\ref{phi21}) we conclude that
$K_{1}=1=K_{2}$ and therefore%
\begin{equation}
\Phi_{5}(x,n)=\Phi_{1}(x,n)+\Phi_{2}(x,n). \label{phi5}%
\end{equation}

Since $-\sqrt{2n}<x<\sqrt{2n},$ we set%
\begin{equation}
x=\sqrt{2n}\sin\left(  \theta\right)  ,\quad-\frac{\pi}{2}<\theta<\frac{\pi
}{2}. \label{theta}%
\end{equation}
Using (\ref{theta}) in (\ref{sigma}), we have%
\begin{equation}
\sigma=\sqrt{2n}\cos\left(  \theta\right)  \mathrm{i.} \label{sigma1}%
\end{equation}
Replacing (\ref{sigma1}) in (\ref{Phi1}), we get%
\begin{gather}
\exp\left[  \frac{x^{2}-\sigma x-n}{2}+n\ln\left(  \sigma+x\right)  \right]
=\label{T1}\\
\exp\left\{  \frac{n}{2}\left[  \ln\left(  2n\right)  -\cos\left(
2\theta\right)  \right]  -n\left[  \frac{1}{2}\sin\left(  2\theta\right)
+\theta-\frac{\pi}{2}\right]  \mathrm{i}\right\} \nonumber
\end{gather}
and%
\begin{equation}
\sqrt{\frac{1}{2}\left(  1+\frac{x}{\sigma}\right)  }=\frac{\exp\left(
-\frac{\theta}{2}\mathrm{i}\right)  }{\sqrt{2\cos\left(  \theta\right)  }}.
\label{T2}%
\end{equation}
Similarly, replacing (\ref{sigma1}) in (\ref{Phi2}), we obtain%
\begin{gather}
\left(  -1\right)  ^{n}\exp\left[  \frac{x^{2}+\sigma x-n}{2}+n\ln\left(
\sigma-x\right)  \right]  =\label{T3}\\
\exp\left\{  \frac{n}{2}\left[  \ln\left(  2n\right)  -\cos\left(
2\theta\right)  \right]  +n\left[  \frac{1}{2}\sin\left(  2\theta\right)
+\theta-\frac{\pi}{2}\right]  \mathrm{i}\right\} \nonumber
\end{gather}
and%
\begin{equation}
\sqrt{\frac{1}{2}\left(  1-\frac{x}{\sigma}\right)  }=\frac{\exp\left(
\frac{\theta}{2}\mathrm{i}\right)  }{\sqrt{2\cos\left(  \theta\right)  }}.
\label{T4}%
\end{equation}
Using (\ref{T1})--(\ref{T4}) in (\ref{phi5}), we have%
\begin{gather}
\Phi_{5}\left[  \sqrt{2n}\sin\left(  \theta\right)  ,n\right]  =\sqrt{\frac
{2}{\cos\left(  \theta\right)  }}\exp\left\{  \frac{n}{2}\left[  \ln\left(
2n\right)  -\cos\left(  2\theta\right)  \right]  \right\} \label{Phi51}\\
\times\cos\left\{  n\left[  \frac{1}{2}\sin\left(  2\theta\right)
+\theta-\frac{\pi}{2}\right]  +\frac{\theta}{2}\right\}  ,\quad-\frac{\pi}%
{2}<\theta<\frac{\pi}{2}.\nonumber
\end{gather}

Thus, we have proved the following:

\begin{theorem}
In the region $\left\vert x\right\vert <\sqrt{2n},$ the Hermite polynomials
have the asymptotic representation%
\begin{gather}
H_{n}\left[  \sqrt{2n}\sin\left(  \theta\right)  \right]  \sim\sqrt{\frac
{2}{\cos\left(  \theta\right)  }}\exp\left\{  \frac{n}{2}\left[  \ln\left(
2n\right)  -\cos\left(  2\theta\right)  \right]  \right\}  \label{Oscilatory}%
\\
\times\cos\left\{  n\left[  \frac{1}{2}\sin\left(  2\theta\right)
+\theta-\frac{\pi}{2}\right]  +\frac{\theta}{2}\right\}  ,\quad n\rightarrow
\infty,\nonumber
\end{gather}
with $-\frac{\pi}{2}<\theta<\frac{\pi}{2}.$
\end{theorem}

\begin{figure}[ptb]
\begin{center}
\rotatebox{270} {\resizebox{!}{5in}{\includegraphics{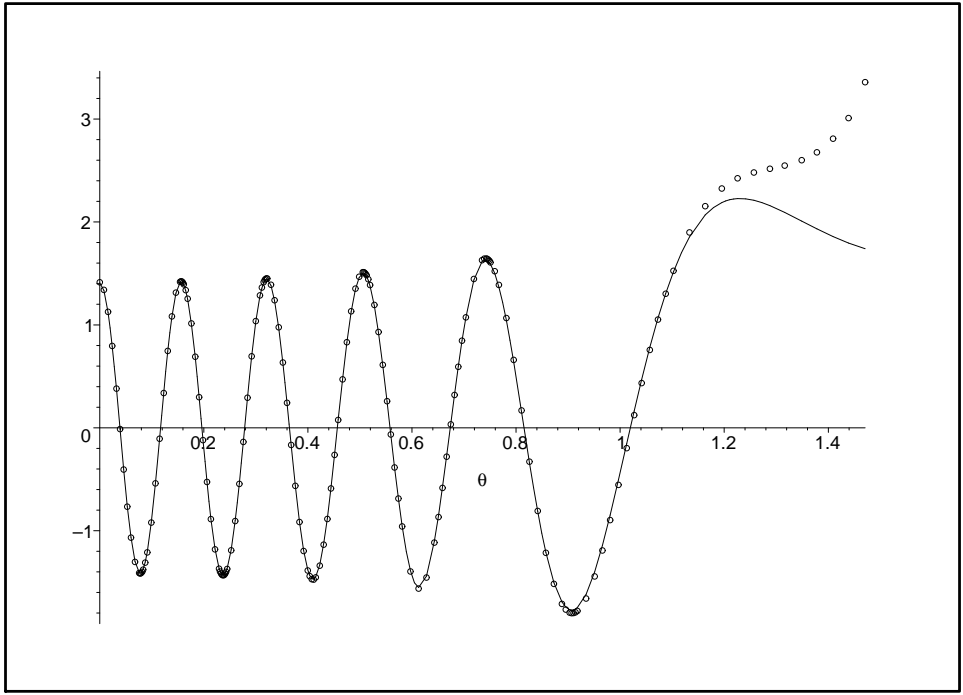}}}
\end{center}
\caption{A comparison of the exact (solid curve) and asymptotic (ooo) values
of $H_{20}\left(  x\right)  $ in the oscillatory region.}%
\label{osci}%
\end{figure}

In Figure \ref{osci} we graph
\[
H_{n}\left[  \sqrt{2n}\sin\left(  \theta\right)  \right]  \exp\left\{
-\frac{n}{2}\left[  \ln\left(  2n\right)  -\cos\left(  2\theta\right)
\right]  \right\}
\]
and
\[
\sqrt{\frac{2}{\cos\left(  \theta\right)  }}\cos\left\{  n\left[  \frac{1}%
{2}\sin\left(  2\theta\right)  +\theta-\frac{\pi}{2}\right]  +\frac{\theta}%
{2}\right\}  ,
\]
with $n=20.$ We only include the range $0\leq\theta<\frac{\pi}{2},$ since both
functions are even.

The same results obtained in this section were derived in \cite{krawherm}
using a different method, based on the limit relation between the Charlier and
Hermite polynomials \cite{koekoek94askeyscheme}.

\section{Zeros}

We shall now find asymptotic formulas for the zeros of the Hermite polynomials
using the results from the previous section. Let's denote by $\zeta_{1}%
^{n}>\zeta_{2}^{n}>\cdots>\zeta_{n}^{n}$ the zeros of $H_{n}\left(  x\right)
,$ enumerated in decreasing order. Then, it follows from (\ref{Oscilatory})
that
\begin{equation}
\zeta_{k}^{n}\sim\sqrt{2n}\sin\left(  \tau_{k}^{n}\right)  ,\quad
n\rightarrow\infty\label{zeta}%
\end{equation}
where $\tau_{k}^{n}$ is a solution of the equation%
\begin{equation}
n\left[  \frac{1}{2}\sin\left(  2\tau_{k}^{n}\right)  +\tau_{k}^{n}-\frac{\pi
}{2}\right]  +\frac{\tau_{k}^{n}}{2}=\left(  1-2k\right)  \frac{\pi}{2}%
,\quad1\leq k\leq n. \label{tau}%
\end{equation}
Solving (\ref{tau}) numerically and using (\ref{zeta}) we get very good
approximations of $\zeta_{k}^{n}$. One could also solve (\ref{tau}) exactly
(as we did in \cite{krawherm}) and obtain a Kapteyn series expansion for
$\tau_{k}^{n}$%
\begin{equation}
\tau_{k}^{n}=\frac{\pi}{2}-\frac{\pi}{2}\left(  4k-1\right)  N^{-1}-%
{\displaystyle\sum\limits_{j=1}^{\infty}}
\frac{1}{j}\mathrm{J}_{j}\left[  \left(  1-N^{-1}\right)  j\right]
\sin\left(  \frac{4k-1}{N}j\pi\right)  , \label{zeros1}%
\end{equation}
where $N=2n+1$ and $\mathrm{J}_{j}\left(  \cdot\right)  $ denotes the Bessel
function of the first kind. However, (\ref{zeros1}) is difficult to analyze
asymptotically. Hence, we will take a different approach and find an
approximation for $\tau_{k}^{n}$ from (\ref{tau}) trough perturbation techniques.

We will consider two cases: $k=O(1)$ which corresponds to the largest zeros of
$H_{n}\left(  x\right)  $ and $k=O\left(  \frac{n}{2}\right)  $, related to
the zeros close to $x=0.$

\subsection{Case I: $k=O(1)$}

Replacing
\begin{equation}
\tau_{k}^{n}=\frac{\pi}{2}-\sum\limits_{i\geq1}^{{}}a_{i}\left(  k\right)
n^{-\frac{i}{3}} \label{tau1}%
\end{equation}
in (\ref{tau}) we obtain, as $n\rightarrow\infty$%
\begin{align}
a_{1}  &  =\frac{1}{2}\kappa^{\frac{1}{3}},\quad a_{2}=-\frac{1}{2}%
\kappa^{-\frac{1}{3}},\quad a_{3}=\frac{\kappa}{120},\quad a_{4}=-\frac
{\kappa^{-\frac{5}{3}}}{30}\left(  \kappa^{2}-5\right) \nonumber\\
a_{5}  &  =\frac{\kappa^{-\frac{7}{3}}}{8400}\left(  3\kappa^{4}+350\kappa
^{2}+1400\right)  ,\quad a_{6}=-\frac{43}{16800}\kappa,\nonumber\\
a_{7}  &  =\frac{\kappa^{-\frac{11}{3}}}{50400}(\kappa^{6}+350\kappa
^{4}-980\kappa^{2}-11200),\nonumber\\
a_{8}  &  =-\frac{\kappa^{-\frac{13}{3}}}{63000}(13\kappa^{6}+475\kappa
^{4}+1400\kappa^{2}+17500),\label{a}\\
a_{9}  &  =\frac{59}{67200}\kappa+\frac{43}{34496000}\kappa^{3},\nonumber\\
a_{10}  &  =-\frac{\kappa^{-\frac{17}{3}}}{1397088000}(23817\kappa
^{8}+2608760\kappa^{6}-4592280\kappa^{4}\nonumber\\
&  -51744000\kappa^{2}-664048000),\nonumber
\end{align}
with
\begin{equation}
\kappa\left(  k\right)  =3\pi\left(  4k-1\right)  . \label{kappa}%
\end{equation}
Using (\ref{tau1})-(\ref{a}) in (\ref{zeta}), we get%
\begin{gather}
\zeta_{k}^{n}\sim\sqrt{2}\left(  n^{\frac{1}{2}}-\frac{\kappa^{\frac{2}{3}}%
}{8}n^{-\frac{1}{6}}+\frac{1}{4}n^{-\frac{1}{2}}-\frac{\kappa^{2}%
+80}{640\kappa^{\frac{2}{3}}}n^{-\frac{5}{6}}\right. \nonumber\\
-\frac{11\kappa^{2}+3920}{179200}n^{-\frac{3}{2}}+\frac{5\kappa^{4}%
+96\kappa^{2}+640}{7680\kappa^{\frac{8}{3}}}n^{-\frac{11}{6}}%
\label{zeroslarge}\\
-\frac{823\kappa^{6}+647200\kappa^{4}-2464000\kappa^{2}-25088000}%
{258048000\kappa^{\frac{10}{3}}}n^{-\frac{13}{6}}\nonumber\\
\left.  +\frac{3064+33\kappa^{2}}{716800}n^{-\frac{5}{2}}\right)  ,\quad
n\rightarrow\infty.\nonumber
\end{gather}

\subsection{Case II: $k=O\left(  \frac{n}{2}\right)  $}

We now set
\begin{equation}
k=\left\lfloor \frac{n}{2}\right\rfloor +1-j=\frac{n}{2}-\alpha+1-j, \label{j}%
\end{equation}
where $\alpha=\operatorname{frac}\left(  \frac{n}{2}\right)  $ (the fractional
part of $\frac{n}{2})$ and $j=0,1,2,\ldots.$ Using (\ref{j}) and%
\begin{equation}
\tau_{k}^{n}=\sum\limits_{i\geq1}^{{}}b_{i}\left(  j\right)  n^{-i}
\label{tau2}%
\end{equation}
in (\ref{tau}) we obtain, as $n\rightarrow\infty$%
\begin{align}
b_{1}  &  =\xi,\quad b_{2}=-\frac{\xi}{4},\quad b_{3}=\frac{\xi}{48}\left(
3+16\xi^{2}\right)  ,\quad\nonumber\\
\quad b_{4}  &  =-\frac{\xi}{192}\left(  3+64\xi^{2}\right)  ,\quad
b_{5}=\frac{\xi}{3840}\left(  15+800\xi^{2}+1024\xi^{4}\right)  ,\label{b}\\
\quad b_{6}  &  =-\frac{\xi}{15360}\left(  15+1600\xi^{2}+7424\xi^{4}\right)
,\nonumber
\end{align}
with%
\begin{equation}
\xi\left(  j\right)  =\frac{\pi}{4}(2j+2\alpha-1). \label{xi}%
\end{equation}
Using (\ref{tau2})-(\ref{b}) in (\ref{zeta}), we obtain%
\begin{equation}
\zeta_{k}^{n}\sim\sqrt{2}\xi\left(  n^{-\frac{1}{2}}-\frac{1}{4}n^{-\frac
{3}{2}}+\frac{3+8\xi^{2}}{48}n^{-\frac{5}{2}}-\frac{3+80\xi^{2}}{192}%
n^{-\frac{7}{2}}\right)  ,\quad n\rightarrow\infty. \label{zeros0}%
\end{equation}

In Table1 we compare the exact value of the positive zeros of $H_{20}\left(
x\right)  $ with the approximations given by solving (\ref{tau}) numerically
and formulas (\ref{zeroslarge}) and (\ref{zeros0}). Note that the biggest
error corresponds to the larger zero, where the asymptotic approximation
(\ref{Oscilatory}) almost breaks down.

\begin{table}[ptb]
\caption{A comparison of the exact and approximate values for the positive
zeros of $H_{20}\left(  x\right)  $.}
\begin{center}
\medskip%
\begin{tabular}
[c]{|c|c|c|c|}\hline
$\zeta_{k}^{n}$ & (\ref{tau}) & (\ref{zeros0}) & (\ref{zeroslarge})\\\hline
.24534 & .24536 & .24536 & \_\\\hline
.73747 & .73751 & .73750 & \_\\\hline
1.2341 & 1.2342 & 1.2340 & \_\\\hline
1.7385 & 1.7387 & 1.7376 & \_\\\hline
2.2550 & 2.2552 & 2.2512 & 2.2592\\\hline
2.7888 & 2.7892 & 2.7779 & 2.7912\\\hline
3.3479 & 3.3486 & \_ & 3.3492\\\hline
3.9448 & 3.9456 & \_ & 3.9460\\\hline
4.6037 & 4.6056 & \_ & 4.6055\\\hline
5.3875 & 5.3939 & \_ & 5.3937\\\hline
\end{tabular}
\end{center}
\end{table}

We summarize the results of this section in the following theorem.

\begin{theorem}
Letting $\zeta_{1}^{n}>\zeta_{2}^{n}>\cdots>\zeta_{n}^{n}$ be the zeros of
$H_{n}\left(  x\right)  ,$ enumerated in decreasing order, we have:

\begin{enumerate}
\item
\[
\zeta_{k}^{n}\sim\sqrt{2}\left(  n^{\frac{1}{2}}-\frac{\kappa^{\frac{2}{3}}%
}{8}n^{-\frac{1}{6}}+\frac{1}{4}n^{-\frac{1}{2}}-\frac{\kappa^{2}%
+80}{640\kappa^{\frac{2}{3}}}n^{-\frac{5}{6}}\right)  ,\quad n\rightarrow
\infty,
\]
where $k=O(1)$ and $\kappa\left(  k\right)  $ was defined in (\ref{kappa}).

\item
\[
\zeta_{k}^{n}\sim\sqrt{2}\xi\left(  n^{-\frac{1}{2}}-\frac{1}{4}n^{-\frac
{3}{2}}+\frac{3+8\xi^{2}}{48}n^{-\frac{5}{2}}-\frac{3+80\xi^{2}}{192}%
n^{-\frac{7}{2}}\right)  ,\quad n\rightarrow\infty,
\]
where $k=$ $\frac{n}{2}-\alpha+1-j,$ $\alpha=\operatorname{frac}\left(
\frac{n}{2}\right)  $ and $\xi\left(  j\right)  $ was defined in (\ref{xi}).
\end{enumerate}
\end{theorem}


\begin{thebibliography}{99}                                                                                               %
\bibitem {MR502800}S.~Ahmed. \newblock Systems of nonlinear equations for the
zeros of {H}ermite  polynomials. \newblock {\em Lett. Nuovo Cimento (2)},
22(9):367--370, 1978.

\bibitem {MR2059744}I.~Area, D.~K. Dimitrov, E.~Godoy, and A.~Ronveaux.
\newblock Zeros of {G}egenbauer and {H}ermite polynomials and connection
coefficients. \newblock {\em Math. Comp.}, 73(248):1937--1951 (electronic), 2004.

\bibitem {MR0457815}F.~Calogero. \newblock On the zeros of {H}ermite
polynomials. \newblock {\em Lett. Nuovo Cimento (2)}, 20(14):489--490, 1977.

\bibitem {MR522774}F.~Calogero and A.~M. Perelomov. \newblock Asymptotic
density of the zeros of {H}ermite polynomials of diverging  order, and related
properties of certain singular integral operators.
\newblock {\em Lett. Nuovo Cimento (2)}, 23(18):650--652, 1978.

\bibitem {MR1373150}O.~Costin and R.~Costin. \newblock Rigorous {WKB} for
finite-order linear recurrence relations with  smooth coefficients.
\newblock {\em SIAM J. Math. Anal.}, 27(1):110--134, 1996.

\bibitem {MR1270235}J.~S. Dehesa, F.~Dom{\'{\i}}nguez-Adame, E.~R. Arriola,
and A.~Zarzo. \newblock Hydrogen atom and orthogonal polynomials. \newblock In
\emph{Orthogonal polynomials and their applications (Erice, 1990)},  volume~9
of \emph{IMACS Ann. Comput. Appl. Math.}, pages 223--229. Baltzer,  Basel, 1991.

\bibitem {MR1342385}H.~Dette and W.~J. Studden. \newblock Some new asymptotic
properties for the zeros of {J}acobi, {L}aguerre,  and {H}ermite polynomials.
\newblock {\em Constr. Approx.}, 11(2):227--238, 1995.

\bibitem {MR0225511}R.~B. Dingle and G.~J. Morgan. \newblock {${\rm WKB}$}
methods for difference equations. {I}, {II}. \newblock {\em Appl. Sci. Res.},
18:221--237; 238--245, 1967/1968.

\bibitem {krawherm}D.~Dominici. \newblock Asymptotic analysis of the
Askey-scheme {II}: from {C}harlier to {H}ermite. \newblock Submitted, 2005.
\newblock arXiv: math.CA/0508264.

\bibitem {MR0361328}T.~Dosdale, G.~Duggan, and G.~J. Morgan.
\newblock Asymptotic solutions to differential-difference equations.
\newblock {\em J. Phys. A}, 7:1017--1026, 1974.

\bibitem {MR2103370}J.~S. Geronimo, O.~Bruno, and W.~Van~Assche.
\newblock W{KB} and turning point theory for second-order difference
equations. \newblock In \emph{Spectral methods for operators of mathematical
physics},  volume 154 of \emph{Oper. Theory Adv. Appl.}, pages 101--138.
Birkh\"auser,  Basel, 2004.

\bibitem {MR1164992}J.~S. Geronimo and D.~T. Smith. \newblock W{KB}
({L}iouville-{G}reen) analysis of second order difference  equations and
applications. \newblock {\em J. Approx. Theory}, 69(3):269--301, 1992.

\bibitem {MR0026415}R.~E. Greenwood and J.~J. Miller. \newblock Zeros of the
{H}ermite polynomials and weights for {G}auss'  mechanical quadrature formula.
\newblock {\em Bull. Amer. Math. Soc.}, 54:765--769, 1948.

\bibitem {MR1465675}K.~Hannabuss.
\newblock {\em An introduction to quantum theory}, volume~1 of \emph{Oxford
Graduate Texts in Mathematics}. \newblock The Clarendon Press Oxford
University Press, New York, 1997.

\bibitem {MR0086169}L.~O. Heflinger. \newblock The asymptotic behaviour of the
{H}ermite polynomials.
\newblock {\em Nederl. Akad. Wetensch. Proc. Ser. A. {\bf 59} = Indag. Math.},
18:255--264, 1956.

\bibitem {koekoek94askeyscheme}R.~Koekoek and R.~F. Swarttouw. \newblock The
{A}skey-scheme of hypergeometric orthogonal polynomials and its  $q$-analogue.
\newblock Technical Report 98-17, Delft University of Technology, 1998. \newblock http://aw.twi.tudelft.nl/~koekoek/askey/.

\bibitem {MR659409}C.~G. Lange and R.~M. Miura. \newblock Singular
perturbation analysis of boundary value problems for  differential-difference
equations. \newblock {\em SIAM J. Appl. Math.}, 42(3):502--531, 1982.

\bibitem {MR804001}C.~G. Lange and R.~M. Miura. \newblock Singular
perturbation analysis of boundary value problems for  differential-difference
equations. {II}. {R}apid oscillations and resonances.
\newblock {\em SIAM J. Appl. Math.}, 45(5):687--707, 1985.

\bibitem {MR804002}C.~G. Lange and R.~M. Miura. \newblock Singular
perturbation analysis of boundary value problems for  differential-difference
equations. {III}. {T}urning point problems.
\newblock {\em SIAM J. Appl. Math.}, 45(5):708--734, 1985.

\bibitem {MR1093425}C.~G. Lange and R.~M. Miura. \newblock Singular
perturbation analysis of boundary value problems for  differential-difference
equations. {IV}. {A} nonlinear example with layer  behavior.
\newblock {\em Stud. Appl. Math.}, 84(3):231--273, 1991.

\bibitem {MR1260916}C.~G. Lange and R.~M. Miura. \newblock Singular
perturbation analysis of boundary value problems for  differential-difference
equations. {V}. {S}mall shifts with layer behavior.
\newblock {\em SIAM J. Appl. Math.}, 54(1):249--272, 1994.

\bibitem {MR1260917}C.~G. Lange and R.~M. Miura. \newblock Singular
perturbation analysis of boundary value problems for  differential-difference
equations. {VI}. {S}mall shifts with rapid  oscillations.
\newblock {\em SIAM J. Appl. Math.}, 54(1):273--283, 1994.

\bibitem {MR0350075}N.~N. Lebedev.
\newblock {\em Special functions and their applications}. \newblock Dover
Publications Inc., New York, 1972.

\bibitem {MR1723073}J.~L. L{\'o}pez and N.~M. Temme. \newblock Hermite
polynomials in asymptotic representations of generalized  {B}ernoulli,
{E}uler, {B}essel, and {B}uchholz polynomials.
\newblock {\em J. Math. Anal. Appl.}, 239(2):457--477, 1999.

\bibitem {MR0109898}F.~W.~J. Olver. \newblock Uniform asymptotic expansions
for {W}eber parabolic cylinder  functions of large orders.
\newblock {\em J. Res. Nat. Bur. Standards Sect. B}, 63B:131--169, 1959.

\bibitem {MR1343543}G.~Pittaluga and L.~Sacripante. \newblock Bounds for the
zeros of {H}ermite polynomials. \newblock {\em Ann. Numer. Math.},
2(1-4):371--379, 1995. \newblock Special functions (Torino, 1993).

\bibitem {MR1509395}M.~Plancherel and W.~Rotach. \newblock Sur les valeurs
asymptotiques des polynomes d'{H}ermite {$H\sb  n(x)=(-I)\sp n
e\sp {\frac{{x\sp 2}}{2}} \frac{{d\sp n }}{{dx\sp n }}\left(
{e\sp {-\frac{{x\sp 2}}{2}}}\right) $}. \newblock {\em Comment. Math. Helv.},
1(1):227--254, 1929.

\bibitem {MR1343062}P.~E. Ricci. \newblock Improving the asymptotics for the
greatest zeros of {H}ermite  polynomials. \newblock {\em Comput. Math. Appl.},
30(3-6):409--416, 1995.

\bibitem {MR0036084}J.~B. Rosser. \newblock Note on zeros of the {H}ermite
polynomials and weights for {G}auss'  mechanical quadrature formula.
\newblock {\em Proc. Amer. Math. Soc.}, 1:388--389, 1950.

\bibitem {MR0048901}H.~E. Salzer, R.~Zucker, and R.~Capuano. \newblock Table
of the zeros and weight factors of the first twenty {H}ermite  polynomials.
\newblock {\em J. Research Nat. Bur. Standards}, 48:111--116, 1952.

\bibitem {MR1303614}P.~R. Subramanian. \newblock Zeros of the {H}ermite
polynomials and the simple {L}aguerre  polynomials are irrational.
\newblock In \emph{Proceedings of the Third Annual Conference of Vijnana
Parishad of India and the National Symposium of Ancient Science in India
(Srinagar, 1993)}, volume~1, pages 49--57, 1993.

\bibitem {MR1329018}P.~R. Subramanian. \newblock Nonzero zeros of the
{H}ermite polynomials are irrational. \newblock {\em Fibonacci Quart.},
33(2):131--134, 1995.

\bibitem {MR1805994}N.~M. Temme and J.~L. L{\'o}pez. \newblock The role of
{H}ermite polynomials in asymptotic analysis. \newblock In \emph{Special
functions (Hong Kong, 1999)}, pages 339--350. World  Sci. Publishing, River
Edge, NJ, 2000.

\bibitem {MR924562}A.~V. Voznyuk. \newblock Asymptotic formulas for parabolic
cylindrical functions and for  {H}ermite polynomials for large values of the
argument. \newblock {\em Vychisl. Prikl. Mat. (Kiev)}, (63):19--24, 125, 1987.

\bibitem {MR1896091}Z.~Wang and R.~Wong. \newblock Uniform asymptotic
expansion of {$J\sb \nu(\nu a)$} via a difference  equation.
\newblock {\em Numer. Math.}, 91(1):147--193, 2002.

\bibitem {MR1971216}Z.~Wang and R.~Wong. \newblock Asymptotic expansions for
second-order linear difference equations  with a turning point.
\newblock {\em Numer. Math.}, 94(1):147--194, 2003.

\bibitem {MR2114641}Z.~Wang and R.~Wong. \newblock Linear difference equations
with transition points. \newblock {\em Math. Comp.}, 74(250):629--653
(electronic), 2005.

\bibitem {MR804827}P.~Wilmott. \newblock A note on the {WKB} method for
difference equations. \newblock {\em IMA J. Appl. Math.}, 34(3):295--302, 1985.

\bibitem {MR0146423}M.~Wyman. \newblock The asymptotic behaviour of the
{H}ermite polynomials. \newblock {\em Canad. J. Math.}, 15:332--349, 1963.
\end{thebibliography}
\end{document}